\documentclass[11pt,a4paper,reqno]{amsart}
\usepackage[T1]{fontenc}
\usepackage{lmodern,microtype,amsmath,amssymb,amsthm,mathtools}
\usepackage[margin=29mm,top=26mm,bottom=28mm]{geometry}
\usepackage[colorlinks=true,allcolors=black,pdfauthor={Zikang Dong and Shengbo Zhao},pdftitle={Large values of quadratic Dirichlet L-functions at bounded height}]{hyperref}
\numberwithin{equation}{section}
\newtheorem{theorem}{Theorem}[section]

\newtheorem{lemma}[theorem]{Lemma}
\newtheorem{proposition}[theorem]{Proposition}
\newcommand{\M}{\mathcal M}
\newcommand{\PP}{\mathcal P}

\newcommand{\dd}{\mathrm d}
\newcommand{\eps}{\varepsilon}

\DeclareMathOperator{\lcm}{lcm}
\newcommand{\half}{\tfrac12}
\allowdisplaybreaks[1]
\title{Large values of quadratic Dirichlet $L$-functions\\ at bounded height}
\author{Zikang Dong}
\author{Weijia Wang}
\author{Hao Zhang}
\author{Shengbo Zhao}
\address[Zikang Dong]{School of Mathematical Sciences, Soochow University,
 Suzhou 215006, P. R. China}

\address[Weijia Wang]{School of Mathematics, Shandong University,
  Jinan, Shandong 250100, P. R. China}
\address[Hao Zhang]{School of Mathematics, Hunan University, Changsha
 410082, P. R. China}
\address[Shengbo Zhao]{School of Mathematical Sciences, Key Laboratory of
 Intelligent Computing and Applications (Tongji University), Ministry of
 Education, Tongji University, Shanghai 200092, China}

\email{zikangdong@gmail.com}
\email{weijiawang@amss.ac.cn}
\email{zhanghaomath@hnu.edu.cn}
\email{shengbozhao@hotmail.com}
\date{}
\begin{document}
\begin{abstract}
Let $\chi_p$ be the quadratic character modulo a prime $p\equiv1\pmod4$. We show
that for any fixed $T>0$, the maximum of $\log|L(\tfrac12+it,\chi_p)|$ over
$X<p\le2X$ and $|t|\le T$ is at least
$(\sqrt{8T^2/(1+16T^2)}+o(1))\sqrt{\log X\,\log_3X/\log_2X}$ as $X\to\infty$. We
use the zero-free half-plane $\Re s>7/8$ for cancellation in a smoothed prime
average. A Gaussian weight keeps the contribution of a squarefree Gál sum, and
a fourth-moment estimate bounds the tail $|t|>T$. We also show that the same
zero-free region implies a sharp cutoff estimate for character sums over
fundamental discriminants, similar to a lemma of Darbar and Maiti.
\end{abstract}
\maketitle

\section{Introduction and results}
Let $\log_j$ be the $j$-fold iterated natural logarithm. Let $X$ be large and
define
\[
 B(X)=\sqrt{\frac{\log X\,\log_3X}{\log_2X}},\qquad
 \PP(X)=\{p\text{ prime}:X<p\le2X,\ p\equiv1\pmod4\}.
\]
For $T\ge0$, define
\[
 M(X,T)=\max_{p\in\PP(X)}\ \max_{|t|\le T}
 |L(\half+it,\chi_p)|.
\]
A prime \(p\in\PP(X)\) is a positive fundamental discriminant. The Kronecker
character \(\chi_p(n)=(\frac{p}{n})\) is the Legendre character
\((\frac{n}{p})\).

\begin{theorem}\label{thm:main}
For any fixed $T>0$,
\begin{equation}\label{eq:main}
 \liminf_{X\to\infty}\frac{\log M(X,T)}{B(X)}
 \ \ge\ C(T),\qquad C(T)=\sqrt{\frac{8T^2}{1+16T^2}}.
\end{equation}
\end{theorem}
In place of \(M(X,T)\) in \eqref{eq:main} one may use
\(\max_{p\in\PP(X),\,0\le t\le T}\Re L(\half+it,\chi_p)\). This maximum is
positive for all sufficiently large \(X\). The formula gives \(C(1/4)=1/2\) and
\(C(1)=\sqrt{8/17}\). For any fixed \(0<c<1/\sqrt2\), the inequality
\(M(X,T)\ge\exp\{cB(X)\}\) holds for all sufficiently large \(X\) if
\[
 T>\frac{c}{\sqrt{8(1-2c^2)}}.
\]

If $H(X)\to\infty$, then for any fixed $T$ we eventually have
$M(X,H(X))\ge M(X,T)$. It follows from \eqref{eq:main}, on letting $T\to\infty$,
that $\liminf_{X\to\infty}\log M(X,H(X))/B(X)\ge1/\sqrt2$. As every $p\in\PP(X)$
is a positive fundamental discriminant, the same lower bounds apply when the
maxima are taken over all fundamental discriminants $X<|d|\le2X$.

At the central point $t=0$, Soundararajan \cite{Sound} proved an unconditional
lower bound of size $\exp\{(1/\sqrt5+o(1))\sqrt{\log X/\log_2X}\}$ for the
maximum over fundamental discriminants. Darbar and Maiti \cite{DM} proved that
the maximum at the central point is at least $\exp\{(1/2+o(1))B(X)\}$, assuming
GRH. Under the same hypothesis, Dong, Wang, Zhang and Zhao
\cite[Theorem~1.1]{DWZZ} improved the constant from $1/2$ to $1$ by using the
square-free Gál sum of de la Bretèche and Tenenbaum \cite{BT}. See
\cite{DWZZprime,FHX} for central values in prime-conductor families. In
Theorem~\ref{thm:main} the ordinate may depend on the character and on $X$. We
discuss why we need $T>0$ in our argument in Section~\ref{sec:height}.

For $1/2\le\theta<1$, let $\mathrm{ZF}(\theta)$ stand for the assertion that
$L(s,\xi)\ne0$ for $\Re s>\theta$ if $\xi$ is a primitive real nonprincipal
Dirichlet character. The same half-plane is then zero-free for the induced
characters because the extra Euler factors are zero only on $\Re s=0$. According
to \cite[Theorem~1.1]{ZFR}, $\mathrm{ZF}(7/8)$ holds. Under the assumption
$\mathrm{ZF}(\theta)$ we show that
\begin{equation}\label{eq:general}
 \liminf_{X\to\infty}\frac{\log M(X,T)}{B(X)}
 \ge 2\sqrt{\frac{16T^2(1-\theta)}{1+16T^2}}
 \qquad (T>0\text{ fixed}).
\end{equation}

The estimate for the sharp discriminant sum in \cite[Lemma~1]{DM} depends on the
square-free and square parts of $n$. If $n=n_0n_1^2$, where $n_0$ is
square-free, then we define
\[
 g_\sigma(n_1)=\sum_{h\mid n_1}\frac{\mu^2(h)}{h^\sigma}
 =\prod_{\ell\mid n_1}(1+\ell^{-\sigma}).
\]
In Euler products, $\ell$ always stands for a prime.
\begin{theorem}\label{thm:sharp}
Let $\mathrm{ZF}(\theta)$ hold and fix $\theta<\sigma<1$. Then there are $C>0$
and $0<\delta<1$, both depending only on $\theta,\sigma$, such that uniformly
for $X\ge2$ and every positive integer $n=n_0n_1^2$,
\begin{align}
 \sum_{\substack{0<|d|\le X\\d\ \mathrm{fundamental}}}\chi_d(n)
 &=\mathbf1_{n_0=1}\frac{X}{\zeta(2)}
     \prod_{\ell\mid n}\frac{\ell}{\ell+1}\notag\\
 &\quad+O_{\theta,\sigma}\!\left(X^\sigma
   \exp\!\big\{C(\log(3n_0))^{1-\delta}\big\}\,
   g_\sigma(n_1)\right).
 \label{eq:sharp}
\end{align}
\end{theorem}
For $\theta=7/8$, we can choose $\sigma=7/8+\eps$ for any fixed $0<\eps<1/8$.
The value of $\delta$ depends on the position of $\sigma$ with respect to the
zero-free boundary. When $n$ is a square, the local densities in the main term
agree with those in the GRH estimate of Darbar and Maiti. In this case the
counting argument gives the remainder $O_\sigma(X^\sigma g_\sigma(n_1))$
unconditionally for every fixed $1/2<\sigma<1$.

\section{Character sums in a zero-free half-plane}\label{sec:characters}
In this section all implicit constants may vary with the fixed lines of
integration.
\begin{lemma}\label{lem:analytic}
Suppose $\mathrm{ZF}(\theta)$ and let $\theta<\sigma<1$ be fixed. Let $\xi$ be a
real nonprincipal character mod $q$. For $\sigma\le u\le2$ we have, uniformly in
$q$ and $t$,
\begin{equation}\label{eq:logder}
 \left|\frac{L'}{L}(u+it,\xi)\right|\ll\log(2q(2+|t|)).
\end{equation}
Also, there are $C,\delta>0$ depending only on $\theta,\sigma$ such that
\begin{equation}\label{eq:partials}
 \left|\sum_{j\le U}\xi(j)\right|
 \ll U^\sigma\exp\{C(\log(3q))^{1-\delta}\}\qquad(U\ge1).
\end{equation}
\end{lemma}
\begin{proof}
Let $s_0=2+it$. We consider the function $F(z)=\log L(s_0+z,\xi)$ for
$|z|<2-\theta$, with the branch from the Euler product at $z=0$. Periodicity and
partial summation on any fixed smaller disc yield $|L(s_0+z,\xi)|\ll q(2+|t|)$,
and the Euler product yields $|F(0)|\ll1$. The Borel--Carathéodory theorem
therefore implies $\max_{|z|\le r}|F(z)|\ll_r\log(2q(2+|t|))$ for any fixed
$r<2-\theta$. If we now take $r>2-\sigma$ and use Cauchy's estimate, we arrive
at \eqref{eq:logder}.

Let $\sigma_0=(\theta+\sigma)/2$ and let $2-\sigma_0<r<2-\theta$. On $|z|\le1/2$
the Euler product shows that $|F(z)|\ll1$. Then Hadamard's three-circles theorem
for the radii $1/2,2-\sigma_0,r$ gives
\[
 |L(u+it,\xi)|\le
 \exp\{C(\log(2q(2+|t|)))^{1-\delta}\}
 \qquad (\sigma_0\le u\le2)
\]
for some $0<\delta<1$. By subadditivity of the power $1-\delta$ we obtain, for
any fixed $\kappa>0$,
\[
 |L(u+it,\xi)|\ll_\kappa \exp\{C(\log(3q))^{1-\delta}\}(2+|t|)^\kappa .
\]

For $U\ge3$ let $x=\lfloor U\rfloor+1/2$, $b=1+1/\log x$ and $H=x^2$. Using
truncated Perron inversion, we find
\begin{equation*}
 \sum_{j\le U}\xi(j)=\frac{1}{2\pi i}\int_{b-iH}^{b+iH}\frac{L(s,\xi)x^s}{s} \dd s
 +O(x\log x/H).
\end{equation*}
The error terms here are $O(x/(H|x-j|))$ for $x/2<j<2x$, and for the other $j$
we use $\sum j^{-b}\ll\log x$. We shift the contour to $\Re s=\sigma_0$ and pull
out $\exp\{C(\log(3q))^{1-\delta}\}$. The vertical part gives
$O_\kappa(x^{\sigma_0}H^\kappa\log(2H))$, and the horizontal parts give
$O_\kappa(x^bH^{\kappa-1})$. Both are $O(x^\sigma)$ as long as
$0<\kappa<(\sigma-\sigma_0)/4$. For $1\le U<3$ the trivial estimate suffices.
\end{proof}
For every positive integer $v$, define $\xi_v\bmod 8v$ to be the real character
that equals the Kronecker symbol $(\frac{u}{v})$ on odd $u$ and is zero on even
$u$. When $v$ is not a square, both $\xi_v$ and $\xi_v\chi_4$ are nonprincipal.
Here $\chi_4$ is the nonprincipal character modulo $4$. Indeed, if $v$ has an
odd prime divisor of odd multiplicity, take $u\equiv1\pmod8$ to be a nonresidue
modulo that prime and $1$ modulo each of the other odd prime divisors. If $2$ is
the only prime of odd multiplicity, take $u\equiv5\pmod8$ and $u\equiv1$ modulo
the odd prime divisors. In both cases $u\equiv1\pmod4$ and $\xi_v(u)=-1$.

\begin{proof}[Proof of Theorem~\ref{thm:sharp}]
Let us first assume that $n_0>1$. Each fundamental discriminant has exactly one
of the three shapes $e,4e,8e$, with $e$ odd and square-free of either sign. In
the three cases the conditions on $e$ are $e\equiv1\pmod4$, $e\equiv3\pmod4$ and
no further congruence condition, respectively. Consider the signs separately,
set $m=|e|$ and expand the residue-class indicators in characters modulo $4$. By
multiplicativity of the Kronecker symbol, the left-hand side of \eqref{eq:sharp}
is a bounded linear combination of sums of the form
\[
 \sum_{m\le Y}\mu^2(m)\xi(m)\mathbf1_{(m,n_1)=1},\qquad
 Y\in\{X,X/4,X/8\},\quad \xi\in\{\xi_{n_0},\xi_{n_0}\chi_4\}.
\]
The signs and the factors $4,8$ yield constants of modulus at most one. The
modulus of both characters is $8n_0$. Using Möbius inversion on $(j,n_1)=1$ and
\eqref{eq:partials}, we obtain
\[
 \left|\sum_{j\le U}\xi(j)\mathbf1_{(j,n_1)=1}\right|
 \ll U^\sigma g_\sigma(n_1)\exp\{C(\log(3n_0))^{1-\delta}\}.
\]
Substituting $\mu^2(m)=\sum_{b^2\mid m}\mu(b)$, we use the estimate with
$U=Y/b^2$ for the inner sum over $j=m/b^2$. The sum over $b$ contributes at most
$\zeta(2\sigma)$.

Let $n_0=1$ and set $D_n(Y)=\sum_{m\le Y,(m,2n)=1}\mu^2(m)$. Then for each
positive odd square-free $m$ exactly one of the two signs is $1\pmod4$ and
exactly one is $3\pmod4$. The required count is therefore
\[
 D_n(X)+\mathbf1_{2\nmid n}\{D_n(X/4)+2D_n(X/8)\}.
\]
The Dirichlet series of the summands in $D_n$ is $\zeta(s)H_n(s)$, where
\[
 H_n(s)=\prod_{\ell\mid 2n}(1-\ell^{-s})
        \prod_{\ell\nmid 2n}(1-\ell^{-2s})
       =\sum_{j\ge1}h_n(j)j^{-s}.
\]
We have $\sum_j|h_n(j)|j^{-\sigma}\le\zeta(2\sigma)g_\sigma(2n)$. By
$\lfloor Y/j\rfloor=Y/j+O((Y/j)^\sigma)$ for all $j\ge1$, we obtain
\[
 D_n(Y)=\frac{Y}{\zeta(2)}\prod_{\ell\mid 2n}\frac{\ell}{\ell+1}
        +O_\sigma(Y^\sigma g_\sigma(2n)).
\]
For odd $n$, the factor $1+1/4+2/8$ cancels the local factor $2/3$. For even
$n$, the prime $2$ already divides $n$. As $n=n_1^2$,
$g_\sigma(2n)\le2g_\sigma(n_1)$, which gives the required bound.
\end{proof}

Let the nonzero, nonnegative $\Phi\in C_c^\infty((1,2))$ be fixed and set
\[
 w_p=(\log p)\Phi(p/X),\quad A_X(v)=\sum_{p\in\PP(X)}w_p\chi_p(v),\quad W_X=A_X(1).
\]
Then, by the prime number theorem for progressions mod $4$, we have
$W_X\sim (X/2)\int_1^2\Phi(u) \dd u$.
\begin{lemma}\label{lem:prime}
Let $\mathrm{ZF}(\theta)$ hold. Then for $\theta<\sigma<1$ and any positive
nonsquare $v$,
\[
 |A_X(v)|\ll_{\Phi,\theta,\sigma} X^\sigma\log(2vX).
\]
For square $v$, $0\le A_X(v)\le W_X$, and $A_X(v)=W_X$ if all prime factors of
$v$ are bounded by $X$.
\end{lemma}
\begin{proof}
If $v$ is not a square, expanding the indicator of $p\equiv1\pmod4$ gives us the
two nonprincipal characters $\xi_v$ and $\xi_v\chi_4$. The Mellin inversion
formula applied to either character $\xi$ of modulus $q=8v$ gives
\[
 \sum_{j\ge1}\Lambda(j)\xi(j)\Phi(j/X)
 =\frac{1}{2\pi i}\int_{(2)}-\frac{L'}{L}(s,\xi)\,X^s\widehat\Phi(s) \dd s,
 \qquad \widehat\Phi(s)=\int_0^\infty \Phi(u)u^{s-1} \dd u .
\]
By shifting the contour to $\Re s=\sigma$ and using \eqref{eq:logder} and the
rapid decay of $\widehat\Phi$, we get $O(X^\sigma\log(2q))$. The error from
subtracting the contribution of higher prime powers is $O_\Phi(X^{1/2}\log^2X)$.
If $v$ is square, then $\chi_p(v)$ is $1$ if $p\nmid v$ and $0$ otherwise.
\end{proof}

\section{The Gaussian weight and the resonating moment}
Let $a,\lambda>0$ be fixed. We write $N=\lfloor X^a\rfloor$ and
$V=\lambda\log X$. For square-free $m,n$, $Q(m,n)=mn/(m,n)^2$ is the smallest
positive $k$ such that $kmn$ is a square. We apply the square-free Gál
construction \cite[Eq.~(1.5) and Section~2]{BT} in the form given in
\cite[Lemma~2.3]{DWZZ}. This gives a set $\M$ of square-free integers. Its
cardinality satisfies $1\le D\le N$, and
\begin{equation}\label{eq:gal}
 y:=\max_{m\in\M}P^+(m)\le (\log N)^{1+o(1)},\qquad
 \sum_{m,n\in\M}Q(m,n)^{-1/2}\ge D\exp\{(2+o(1))B(N)\},
\end{equation}
where $P^+(m)$ is the largest prime factor of $m$, with the convention
$P^+(1)=1$. The upper bound on $y$ comes from the prime intervals in the
construction. Thus $\log m\le\sum_{\ell\le y}\log\ell\ll y$ for every $m\in\M$,
where the sum is over primes. We write
\[
 b_V(k)=\exp\!\left(-\frac{(\log k)^2}{4V}\right),\qquad
 P_p(V)=\sum_{k\ge1}\frac{b_V(k)\chi_p(k)}{\sqrt k},\qquad
 R_p=\sum_{m\in\M}\chi_p(m).
\]
The series for $P_p(V)$ converges absolutely for each $V>0$.

\begin{lemma}\label{lem:gal}
With these choices we have, as $X\to\infty$,
\[
 \frac{1}{D}\sum_{m,n\in\M}\frac{b_V(Q(m,n))}{\sqrt{Q(m,n)}}
 \ge \exp\{(2\sqrt a+o(1))B(X)\}.
\]
\end{lemma}
\begin{proof}
Put $Z=\exp((\log X)^{7/10})$ and let $m$ be fixed. The primes dividing $Q(m,n)$
form the symmetric difference of the prime supports of $m$ and $n$ and so
uniquely determine $n$. Thus
\[
 \sum_{n\in\M}Q(m,n)^{-1/3}
 \le\prod_{\ell\le y}(1+\ell^{-1/3})\le \exp\{O(y^{2/3})\}.
\]
Since $2/3<7/10$, the contribution to the unweighted Gál sum from pairs for
which $Q>Z$ is bounded by $D Z^{-1/6}\exp\{O(y^{2/3})\}=o(D)$. For the other
pairs we have
\[
 -\log b_V(Q)\le\frac{(\log Z)^2}{4V}
 =\frac{(\log X)^{2/5}}{4\lambda}=o(B(X)).
\]
The lower bound follows from \eqref{eq:gal} and $B(N)=(\sqrt a+o(1))B(X)$.
\end{proof}

\begin{proposition}\label{prop:moment}
Suppose $\mathrm{ZF}(\theta)$ and
\begin{equation}\label{eq:budget}
 a+\frac{\lambda}{4}<1-\theta.
\end{equation}
Then we have $S_1:=\sum_{p\in\PP(X)} w_pR_p^2\sim W_XD$ and
\begin{equation}\label{eq:ratio}
 \frac{1}{S_1}\sum_{p\in\PP(X)}w_pR_p^2 P_p(V)
 \ge \exp\{(2\sqrt a+o(1))B(X)\}.
\end{equation}
\end{proposition}
\begin{proof}
In $S_1$, $mn$ is a square if and only if $m=n$. In this case all its prime
factors are at most $y<X$, and each diagonal entry contributes $W_X$. As
$\log(2Xmn)\ll\log X+y$, Lemma~\ref{lem:prime} gives, for every sufficiently
small fixed $\eta>0$,
\[
 S_1=W_XD+O_\eta(X^{\theta+\eta}D^2).
\]
For nonsquare $kmn$ in the numerator we use Lemma~\ref{lem:prime} with
$\sigma=\theta+\eta/2$. We split the $k$-sum into the intervals
$e^j\le k<e^{j+1}$ for $j\ge0$. The weighted sum of $\log(2Xkmn)$ is bounded by
a constant times
\[
 \sum_{j\ge0}(\log(2Xmn)+j+1)\exp\!\left(\frac j2-\frac{j^2}{4V}\right)
 \ll e^{V/4}(1+V)^{1/2}(\log(2Xmn)+V+1).
\]
The last bound is obtained by completing the square at $j=V$. By
$\log(2Xmn)\ll\log X+y$, the factors other than $e^{V/4}$ are
$O_\eta(X^{\eta/2})$ with an implied constant independent of $m,n$. Hence the
total absolute contribution of nonsquare $kmn$ is
\begin{equation*}
 \ll_\eta X^{\theta+\eta}D^2 e^{V/4}.
\end{equation*}
As $D\le X^a$ and $e^{V/4}=X^{\lambda/4}$, it follows from \eqref{eq:budget}
that both errors are $O(X^{1-\delta_1}D)$ for some $\delta_1>0$ when $\eta$ is
sufficiently small.

The square terms are nonnegative. For each pair $m,n$ we keep the term
$k=Q(m,n)$. Then $kmn=\lcm(m,n)^2$ has all prime divisors no larger than $y$ and
has prime average $W_X$. Dividing by $S_1\sim W_XD\asymp XD$ and using
Lemma~\ref{lem:gal} gives us \eqref{eq:ratio}.
\end{proof}

\section{Restriction to a fixed interval of heights}\label{sec:height}
For the Gaussian weight $b_V$ we use the Fourier transform of a probability
density. Indeed,
\[
 G_V(t)=\sqrt{\frac V\pi}\,e^{-Vt^2},\qquad
 \int_{\mathbb R}G_V(t)e^{-it\log k} \dd t=b_V(k).
\]
Consequently
\begin{equation}\label{eq:identity}
 P_p(V)=\int_{\mathbb R}G_V(t)L(\half+it,\chi_p) \dd t .
\end{equation}
The partial sums $\sum_{k\le K}\chi_p(k)k^{-1/2-it}$ converge pointwise to
$L(\half+it,\chi_p)$. By partial summation, the estimate
$|\sum_{k\le u}\chi_p(k)|\le p$ gives an $O(p(1+|t|))$-bound for these sums,
uniform in $K$. We can therefore take the Gaussian transform of the partial sums
and pass to the limit by dominated convergence.

We control the contribution of $|t|>T$ to the resonating average by taking
$\pi_p=w_p R_p^2/S_1$, so $\sum_p\pi_p=1$. If we restrict Heath-Brown's
fourth-moment estimate for primitive real characters \cite[Theorem~2]{HB} to the
characters $\chi_p$ with $p\in\PP(X)$, then for any real $t$,
\begin{equation}\label{eq:HB}
 \sum_{p\in\PP(X)} |L(\half+it,\chi_p)|^4
 \ll_\eps X^{1+\eps}(2+|t|)^{1+\eps}.
\end{equation}
From $|R_p|\le D$, $w_p\ll\log X$ and $S_1\asymp XD$ we get
$\max_p\pi_p\ll D(\log X)/X$. Applying Hölder's inequality together with
\eqref{eq:HB} and $D\le X^a$ gives, after adjusting $\eps$,
\begin{align}
 \sum_p\pi_p |L(\half+it,\chi_p)|
 &\le\left(\sum_p\pi_p |L(\half+it,\chi_p)|^4\right)^{1/4}\notag\\
 &\ll_\eps X^{a/4+\eps}(2+|t|).
 \label{eq:weighted}
\end{align}
Let $T>0$ be fixed. Then for $V\ge1$ we have
$\int_{|t|>T}G_V(t)(2+|t|) \dd t\ll_T V^{-1/2}e^{-VT^2}$. Integrating
\eqref{eq:weighted} against $G_V$ over $|t|>T$, we control the tail of the
average \eqref{eq:ratio} in absolute value by
\begin{equation}\label{eq:tail}
 O_{T,\eps}\!\left(X^{a/4+\eps-\lambda T^2}\right).
\end{equation}

\begin{proof}[Proof of \eqref{eq:general} and Theorem~\ref{thm:main}]
Fix $0<a<16T^2(1-\theta)/(1+16T^2)$ and $\lambda$ so that
$a/(4T^2)<\lambda<4(1-\theta-a)$. Then \eqref{eq:budget} is satisfied and the
tail \eqref{eq:tail} is $o(1)$ for small enough $\eps>0$. Now we average
\eqref{eq:identity} with the nonnegative weights $\pi_p$ and take real parts.
Because $\Re L(\half+it,\chi_p)$ is even in $t$, $\sum_p\pi_p=1$ and
$\int_{-T}^T G_V(t) \dd t\le1$, the integral on $[-T,T]$ is bounded above by
\[
 \max\!\left\{0,\max_{p\in\PP(X),\,0\le t\le T}\Re L(\half+it,\chi_p)\right\}.
\]
Using \eqref{eq:ratio} and \eqref{eq:tail}, we see that this maximum is bounded
below by $\exp\{(2\sqrt a+o(1))B(X)\}$. So the inner maximum is eventually
positive and the same bound applies to $M(X,T)$ as well. Passing to the limit
inferior and letting $a$ tend to its upper bound, we get \eqref{eq:general}.
Taking $\theta=7/8$ from \cite[Theorem~1.1]{ZFR} gives \eqref{eq:main}.
\end{proof}

The tail condition $\lambda T^2>a/4$ implies $T>0$ for $a>0$. The approximate
functional equation at $t=0$ for conductors of size $X$ has length
$X^{1/2+o(1)}$ and sum of absolute coefficients $O(X^{1/4+\eps})$. The error
term from applying Lemma~\ref{lem:prime} term by term to the nonsquares in the
resonating first moment (see \cite[Section~3]{DWZZ}) is
$O_\eps(X^{\theta+1/4+\eps}D^2)$. We can use Mellin inversion for the
conductor-dependent weight, or use partial summation with a smooth prime cutoff
and uniform derivative bounds for the weight. For $\theta=7/8$ this error bound
is larger than $XD$ even for $D=1$. So using this method to reach the scale of
Theorem~\ref{thm:main} at $t=0$ would require further cancellation in the
central-value moment.

\bigskip
\subsection*{Declaration of AI Use}
Large language models were used to assist in the preparation of
this manuscript, including exploration
of mathematical ideas and language polishing. The authors have reviewed and verified
all AI-assisted content and take full responsibility for the
accuracy, originality, and mathematical correctness of the work.

\end{document}